\documentclass[11pt]{amsart}
\usepackage[top=1in, bottom=1in, left=1in, right=1in]{geometry}
\usepackage{mathrsfs}
\usepackage{amsfonts}
\usepackage{amsmath}
\usepackage{amssymb}
\usepackage{bbm}
\usepackage{physics}
\numberwithin{equation}{section}
\usepackage{times}
\usepackage[usenames,dvipsnames]{color}
\usepackage{comment}
\usepackage{mathtools}
\usepackage{bbm}
\usepackage{ dsfont }
\usepackage{bm}
\usepackage{esvect}
\usepackage{cite}
\usepackage{hyperref}
\hypersetup{
colorlinks = true,
linkcolor={black},
urlcolor={blue},
citecolor={blue},    
urlcolor = {blue},
citebordercolor = {0.33 .58 0.33},
 linkbordercolor = {0.99 .28 0.23},
 breaklinks=true}
 
 \newcommand{\Mod}[1]{\ (\mathrm{mod}\ #1)}

\newcommand{\kommentar}[1]{}

\newtheorem{thm}{Theorem}[section]

\newtheorem{lem}[thm]{Lemma}

\theoremstyle{remark}
\newtheorem{remark}{Remark}[section]

\newcommand{\Li}{\mathrm{Li}} 

\usepackage{geometry}
\title[Least Totients]
{Smallest totient in a residue class}

\author{Abhishek Jha}
\address{
Department of Mathematics,
University of Illinois Urbana-Champaign,
Altgeld Hall, 1409 W. Green Street,
Urbana, IL, 61801, USA}
\email{jha33@illinois.edu}

\DeclareMathOperator*{\f}{{f}}

\begin{document}

\maketitle

\begin{abstract}
We obtain a totient analogue for Linnik's theorem in arithmetic progressions. Specifically, for any coprime pair of positive integers $(m,a)$ such that $m$ is odd, there exists $n\le m^{2+o(1)}$ such that  $\varphi(n)\equiv a\Mod{m}$. 
\end{abstract}

\section{Introduction}
Let $\varphi(n)$ be the Euler function of $n$. A totient is defined as an integer that appears as a value of $\varphi(n)$. The literature is rich with results exploring the distribution of totients in residue classes. Dence and Pomerance \cite{MR1642868} established that if a congruence class $a\Mod m$ contains at least one multiple of 4, it contains infinitely many totients. Later, Ford, Konyagin, and Pomerance \cite{MR1689545} showed that almost all even integers which are $2\Mod{4}$ lie in a residue class that is totient-free.

Let $m$ and $a$ be two coprime positive integers such that $m$ is odd. Let  $N(a,m)$ denote the smallest positive integer $n$ for which $\varphi(n)\equiv a\Mod{m}$. The problem of upper-bound estimates for $N(a, m)$ has garnered significant attention in recent years. These results are closely related to the renowned Linnik problem of bounding the least prime $P(a,m)$ in the arithmetic progression $a\Mod{m}$. {By a result of Xylouris \cite{MR2825574}, we know that $P(a,m)\ll m^{5.18}$ when $\gcd(a,m)=1$. Consequently, this implies that $N(a,m)\ll m^{5.18}$ whenever $\gcd(a+1,m)=1$. Under Generalised Riemann Hypothesis, Lamzouri, Li, and Soundararajan have shown in \cite[Corollary 1.2]{MR3356031} that  $P(a,m)\le (\varphi(m)\log m)^{2}$. This leads to the bound $N(a,m)\ll m^{2+\epsilon}$ assuming that $\gcd(a+1,m)=1$.}

It has been shown that with simpler methods, better unconditional upper bounds for $N(a,m)$ can be achieved. The first result in this direction was by Friedlander and Shparlinski \cite{MR2331570} (See also \cite{MR2418809}) who had proven that if $m$ is a prime number, then $N(a,m)\ll_{\epsilon} m^{2.5+\epsilon}$, a result later refined by Garaev \cite{MR2506368} to $N(a,m)\ll_{\epsilon} m^{2+\epsilon}$. In the case when $m$ is composite,  Friedlander and Shparlinski established that for some $A =
A(\epsilon) > 0$ if $\gcd(a, m) = 1$ and if $m$ has no prime divisors $p < (\log m)^{A(\epsilon)}$, then $N(a,m)\ll_{\epsilon} m^{3+\epsilon}$. This was ultimately improved by Cilleruelo and Garaev \cite{Cilleruelo2009} to $N(a,m)\ll_{\epsilon} m^{2+\epsilon}$ for the same set of positive integers $m$.

The goal of the present paper is to obtain the same upper bound uniformly for all odd positive integers $m$. {Note that when $m$ is even, then there exists a totient in reduced residue class $a\Mod{m}$ if and only if $a\equiv1\Mod{m}$. Moreover, there is only one totient $\varphi(1)=1$ in such a residue class.}
\begin{thm}
    For any $\epsilon>0$ and for all odd positive integers $m$ and integer $a$ with $\gcd(a,m)=1$, we have the bound \[N(a,m)\ll_{\epsilon} m^{2+\epsilon}.\]
\end{thm}

\noindent The restriction that $\gcd(a,m)=1$ is crucial, as Friedlander and Luca \cite{MR2471950} have shown that there exists a sequence of arithmetic progressions $a_k\Mod{m_k}$ with $m_k\rightarrow\infty$ as
$k\rightarrow\infty$ and {$\gcd(a_k,m_k)>1$} such that $N(a_k,m_k)$ exists and 
\[\frac{\log N(a_k,m_k)}{\log m_k}\rightarrow\infty\text{ as }k\rightarrow\infty.\]

\noindent {The main arithmetic input in our work is obtaining asymptotic formulas for character sums over shifted primes where the characters have small conductors. Combining these estimates with the ideas of Cillerurlo and Garaev \cite{Cilleruelo2009}, we are able to get a uniform upper bound for $N(a,m)$.}
\subsection{Notations}When there is no danger of confusion, we write $(a,b)$ instead of $\gcd(a,b)$. We write $\tau(n)$ as the number of divisors of $n$. Throughout this paper, the letter $p$ will always denote a prime. {We set $\f(\chi)$ as the conductor of a Dirichlet character $\chi\Mod{m}$. Let $\mathrm{Li}(x)$ denote the logarithmic integral, that is, \[\mathrm{Li}(x)=\int_{2}^{x}\frac{\mathrm{d}t}{\log t}.\]} We employ the Landau–Bachmann “Big Oh” and “little oh” notations $\mathcal{O}$ and o,
as well as the associated Vinogradov symbols $\ll$ and $\gg$, with their usual meanings. Any dependence of the implied constants is indicated with subscripts.

\subsection{Acknowledgements} The author thanks Kevin Ford for suggesting this problem, for helpful discussions, and for careful reading of the manuscript. The author is grateful to Ayan Nath and the anonymous referee for helpful comments on earlier versions of the manuscript.
\section{Setup}
\noindent In what follows, $\mathds{1}_{a,m}$ denotes the indicator function that  $3\mid m$ and $a\equiv 2\Mod{3}$. {We look for a solution of the congruence in question in the form $n = 4^{\mathds{1}_{a,m}}p_1 p_2 p_3$, where $p_j$ are prime numbers that run through certain disjoint intervals.}

Let $k \geq 10$ be a fixed positive integer. Let $I_1, I_2, I_3$ be sets of primes defined as follows:
\begin{align*}
I_1 &= \{p : 0.5m^{1+1/k} < p \leq m^{1+1/k}, (p - 1, m) = 1\}, \\
I_2 &= \{p : 0.5m < p \leq m, (p - 1, m) = 1\}, \\
I_3 &= \{p : 0.5m^{1/k} < p \leq m^{1/k}, (p - 1, m) = 1\}.
\end{align*}

The sets $I_1, I_2, I_3$ are pairwise disjoint for any sufficiently large integer $m$. The following lemma tells us about the cardinality of these sets. 
\begin{lem}\label{lm1}
    We have for $m\ge 1$ and any $C\ge 2$,
    \begin{equation*}
        |I_{j}|=[\mathrm{Li}(f_{j}(m))-\mathrm{Li}(0.5\,f_{j}(m))]\prod_{p\mid m}\left(1-\frac1{p-1}\right)+\mathcal{O}_{C}(f_j(m)\log^{-C}m),
    \end{equation*} where $f_1(m)$, $f_2(m)$ and $f_3(m)$ are equal to $m^{1+1/k}$, $m$ and $m^{1/k}$ respectively.
\end{lem} 
\begin{proof}
    See \cite[Lemma 4]{MR2331570}.
\end{proof}

\noindent We will prove that if $m$ is a large integer such that $(a, m) = 1$, then the congruence
{\begin{equation}\label{eq1}
\varphi(4^{\mathds{1}_{a,m}}p_1 p_2 p_3)=(1+\mathds{1}_{a,m})\,(p_1 - 1)\,(p_2 - 1)\,(p_3 - 1) \equiv a\Mod{m}, \quad p_j \in I_j, \, j = 1, 2, 3
\end{equation}}
has solutions. The number $J$ of solutions of this congruence is equal to
\begin{equation}\label{eq16}
J = \frac{1}{\varphi(m)} \sum_{\chi} \sum_{p_1, p_2, p_3} \chi((p_1 - 1)(p_2 - 1)(p_3 - 1)) \overline{\chi}(a)\chi(1+\mathds{1}_{a,m})
\end{equation}
where $\chi$ runs through all Dirichlet characters modulo $m$ and the primes $p_1, p_2, p_3$ run through the sets $I_1, I_2, I_3$ respectively. If $3\mid m$, let $\psi$ denote the unique character$\Mod{m}$ induced by the nontrivial character$\Mod{3}$. Denote $\mathds{1}_{3\mid m}$ as the indicator function that $3$ divides $m$. Thus
\begin{align}\label{eq36}
J &= \frac{|I_1||I_2||I_3|}{\varphi(m)} +\frac{\mathds{1}_{3\mid m}S_1(\psi)S_2(\psi)S_3(\psi)\overline{\psi}(a)\psi(1+\mathds{1}_{a,m})}{\varphi(m)}\\&+\frac{\theta}{\varphi(m)} \sum_{\chi \neq \chi_0,\psi} |S_1(\chi)||S_2(\chi)||S_3(\chi)|, \quad |\theta| \leq 1,
\end{align}
where
\[
S_j(\chi) = \sum_{p \in I_j} \chi(p - 1), \quad j = 1, 2, 3.
\]
We notice that 
\[ \mathds{1}_{3\mid m}S_1(\psi)S_2(\psi)S_3(\psi)\overline{\psi}(a)\psi(1+\mathds{1}_{a,m})=\begin{cases} 
      |I_1||I_2||I_3| & \text{if }3\mid m, \\
      0 & \text{otherwise.}
   \end{cases}
\]

\noindent Henceforth, it follows that 
\begin{align*}
    J=\frac{(1+\mathds{1}_{3\mid m})|I_1||I_2||I_3|}{\varphi(m)}+\frac{\theta}{\varphi(m)} \sum_{\chi \neq \chi_0,\psi} |S_1(\chi)||S_2(\chi)||S_3(\chi)|, \quad |\theta| \leq 1.
\end{align*}
To prove that $J > 0$, it is enough to prove that
\begin{equation*}
S=\sum_{\chi \neq \chi_0,\psi} |S_1(\chi)||S_2(\chi)||S_3(\chi)| < |I_1||I_2||I_3|.
\end{equation*}

\noindent Let $A=4(k+3)^{2}$. We split the left sum into two subsums as follows: 
\begin{equation}\label{eq3}
    \sum_{\chi \neq \chi_0,\psi} |S_1(\chi)||S_2(\chi)||S_3(\chi)|=S_{\,\f(\chi)\,\le\, \log^{A} m}+S_{\,\f(\chi)\,>\, \log^{A} m},
\end{equation}
where the first sum is taken over characters with small conductors. We deal with the sums separately in upcoming sections.

\section{ Sum involving characters of small conductor}
To evaluate $S_{\,\f(\chi)\le\log^{A} m}$, the first step is to estimate $S_j(\chi)$. Let $\chi_d$ be the primtive character $\mod{d}$ inducing $\chi$ and $\chi_{0,d}$ be the trivial character$\Mod{d}$. We take $m=dr$. Denote $(r,d)_{\infty}:={r}/{\prod_{p\mid (d,r)}p^{\nu_p(r)}}$ {where $\nu_{p}(r)$ is the largest exponent such that $p^{\nu_{p}(r)}$ divides $r$.} Then
\begin{align}\label{eq5}
    S_j(\chi)=\sum_{p\in I_j}\chi(p-1)=\sum_{v \Mod{d}}\chi_{0,d}(v)\,\chi_{d}(v-1)\sum_{\substack{ p\in I_j \\ p\equiv v \Mod{d}\\ (p-1,(r,d)_{\infty})=1}}1
\end{align}
We detect the coprimality condition in the inner sum using M\"obius function, obtaining 
\begin{align*}
    \sum_{\substack{ p\in I_j \\ p\equiv v \Mod{d}\\ (p-1,(r,d)_{\infty})=1}}1&=\sum_{\substack{ p\in I_j \\ p\equiv v \Mod{d}}}\sum_{\substack{c\mid p-1\\ c\mid (r,d)_{\infty}}}\mu(c)=\sum_{c\mid (r,d)_{\infty}}\mu(c)\sum_{\substack{p\in I_j\\ p\equiv v\Mod{d}\\ p\equiv 1\Mod{c}}}1.\\&=\sum_{c\mid (r,d)_{\infty}}\mu(c)\Delta_{j}(m,cd,v_c),
\end{align*}
where { $v_c$ is the unique positive integer such that \begin{equation*}
    (v_c \,\mathrm{mod}\,{dc})=(v\,\mathrm{mod}\,{d})\cap (1\,\mathrm{mod}\,{c}),
\end{equation*} and }$$\Delta_{j}(m,cd,v_c)=\pi(f_{j}(m);dc,v_c)-\pi(0.5f_{j}(m);dc,v_c).$$

\noindent By the trivial bound, $\pi(x;dc,v_c)\le x/dc$, we get that 
\[\sum_{\substack{c\mid (r,d)_{\infty} \\ c\,\ge\, f_{j}(m)^{1/3}/d}}\mu(c)\Delta_{j}(m,dc,v_c)\ll f_{j}(m)^{2/3}\tau(m).\] 
For the remaining sum, we have 
\begin{equation}\label{eq4}
\sum_{\substack{c\mid (r,d)_{\infty} \\ c\,<\, f_{j}(m)^{1/3}/d}}\mu(c)\Delta_{j}(m,dc,v_c)=\frac{\Li(f_{j}(m))-\Li(0.5\,f_{j}(m))}{\varphi(d)}\sum_{\substack{c\mid (r,d)_{\infty}\\ c<f_{j}(m)^{1/3}/d}}\frac{\mu(c)}{\varphi(c)} +\mathcal{O}(R(m)),\end{equation} where the remainder term is given by 
\[R(m)=\sum_{\substack{c\mid (r,d)_{\infty}\\ c<f_{j}(m)^{1/3}/d}}\bigg|\Delta_{j}(m;cd,v_c)-\frac{\mathrm{Li}(f_{j}(m))-\mathrm{Li}(0.5\,f_{j}(m))}{\varphi(cd)}\bigg|\ll_{B} f_{j}(m)\log^{-B} m\] by the Bombieri-Vinogradov theorem for any constant $B$ (see, e.g., \cite[Theorem 18.9]{MR3971232}). Finally, for the main term in $(\ref{eq4})$, we have 
\begin{equation*}
    \sum_{\substack{c\mid (r,d)_{\infty} \\ c<f_{j}(m)^{1/3}/d}}\frac{\mu(c)}{\varphi(c)}=\prod_{p\mid (r,d)_{\infty}}\left(1-\frac1{p-1}\right)+\mathcal{O}\left(\frac{d\tau(m)\log \log m}{f_{j}(m)^{1/3}}\right).
\end{equation*}
We denote $\rho_{\chi_{d}}$ as $$\rho_{\chi_{d}}=\frac{\sum_{v \Mod{d}}\chi_{0,d}(v)\,\chi_{d}(v-1)}{\varphi(d)}.$$
Putting everything in $(\ref{eq5})$, we get 
\begin{align*}
    S_j(\chi)&= \sum_{v\Mod{d}}\chi_{0,d}(v)\chi_{d}(v-1)\left(\frac{\mathrm{Li}(f_{j}(m))-\mathrm{Li}(0.5\,f_{j}(m))}{\varphi(d)}\prod_{p\mid (r,d)_{\infty}}\left(1-\frac1{p-1}\right)+\mathcal{O}_{B}\left(\frac{f_{j}(m)}{\log^{B}m}\right)\right).
\end{align*}
Note that the first term inside brackets is much larger than the second term since $d$ is at most a fixed power of $\log x$. 
This in turn results in the following asymptotic 
\begin{align}  
S_1(\chi)&=|I_1|\prod_{p\mid (r,d)_{\infty}}\left(1-\frac1{p-1}\right)\cdot\prod_{p\mid m}\left(1-\frac{1}{p-1}\right)^{-1}\,\rho_{\chi_{d}}+\mathcal{O}_{B}\left(\frac{m^{1+1/k}}{\log^{B-A}m}\right)\\&=|I_1|\prod_{p\mid d}\left(1+\frac1{p-2}\right)\,\rho_{\chi_{d}}+\mathcal{O}_{B}\left(\frac{m^{1+1/k}}{\log^{B-A}m}\right),\label{eq6}\end{align} by Lemma \ref{lm1}.
Similarly, we can obtain 
\begin{equation}\label{eq8}
     S_2(\chi)= |I_2|\prod_{p\mid d}\left(1+\frac1{p-2}\right)\,\rho_{\chi_{d}}+\mathcal{O}_{B}\left(\frac{m}{\log^{B-A}m}\right),
\end{equation} and 
\begin{equation}\label{eq9}
     S_3(\chi)= |I_3|\cdot\prod_{p\mid d}\left(1+\frac1{p-2}\right)\,\rho_{\chi_{d}}+\mathcal{O}_{B}\left(\frac{m^{1/k}}{\log^{B-A}m}\right).
\end{equation}

Our next goal is to evaluate the constant $\rho_{\chi_d}$. We will follow ideas in \cite[pp.11-12]{pollack2024mean} to obtain the following lemma.
{
\begin{lem}\label{lm}
    Let $\chi_d$ be a primitive Dirichlet character$\Mod{d}$ for some odd positive integer $d>1$. Then, 
    \begin{equation*}
        \rho_{\chi_d}=\mu(d)\chi_{d}(-1)\prod_{p\mid d}(p-1)^{-1}.
    \end{equation*}
\end{lem}
\begin{proof}
We can write $\chi_d$ uniquely in the form $\prod_{p^{\alpha}||d}\chi_{p^\alpha}$ where $\chi_{p^\alpha}$ is a primitive Dirichlet character $\mathrm{mod}\,{p^{\alpha}}$. Therefore, $\rho_{\chi_d}=\prod_{p^{\alpha}||d}\rho_{\chi_{p^{\alpha}}}$, where for each prime power $p^{\alpha}||d$, we have \begin{align}
   \varphi(p^{\alpha})\rho_{\chi_{p^{\alpha}}}&=\sum_{v\Mod{p^{\alpha}}}\chi_{0,p^{\alpha}}(v)\chi_{p^\alpha}(v-1)=\sum_{\substack{v\Mod{p^{\alpha}}\\ (v,p)=1}}\chi_{p^\alpha}(v-1)\\&=\sum_{v\Mod{p^{\alpha}}}\chi_{p^\alpha}(v)-\sum_{\substack{v\Mod{p^{\alpha}}\\ v\equiv -1\Mod{p}}}\chi_{p^\alpha}(v)\label{eq10}.
\end{align}
 The first sum in $(\ref{eq10})$ is just zero. To evaluate the second sum, consider a primitive root $g\Mod{p^{\alpha}}$ which exists since $p$ is odd. We can observe that the residues $v\Mod{p^{\alpha}}$ such that $v\equiv1\Mod{p}$ are a permutation of the residues $\{g^{(p-1)k}\Mod{p^{\alpha}}\,:\,0\le k<p^{\alpha-1}\}$. Therefore,
\begin{align*}
    \sum_{\substack{v\Mod{p^{\alpha}}\\ v\equiv -1\Mod{p}}}\chi_{p^\alpha}(v)&=\chi_{p^\alpha}(-1)\sum_{\substack{v\Mod{p^{\alpha}}\\ v\equiv 1\Mod{p}}}\chi_{p^\alpha}(v)=\chi_{p^\alpha}(-1)\sum_{0\le k< p^{\alpha-1}}\chi_{p^\alpha}(g^{p-1})^{k}\\&=\chi_{p^\alpha}(-1)\mathds{1}_{(\chi_{p^\alpha})^{p-1}=\chi_{0,p^{\alpha}}}\,p^{\alpha-1}=\chi_{p^\alpha}(-1)\mathds{1}_{\f(\chi_{p^\alpha})\mid p}\,p^{\alpha-1}.
\end{align*}
Finally, we obtain that $$ \rho_{\chi_{p^{\alpha}}}={-\chi_{p^{\alpha}}(-1)\mathds{1}_{\f(\chi_{p^{\alpha}})\mid p}}{(p-1)^{-1}},$$ for each prime power $p^{\alpha}||d$. Multiplying these relations over all prime powers, we obtain
\begin{align*}
    \rho_{\chi_d}&=\prod_{p^{\alpha}||d}\rho_{\chi_{p^{\alpha}}}=\mathds{1}_{\mu^{2}(d)=1}\prod_{p^{\alpha}||d}{-\chi_{p^{\alpha}}(-1)}{(p-1)^{-1}}\\&=\mu(d)\chi_{d}(-1)\prod_{p\mid d}(p-1)^{-1},
\end{align*} which is what we wanted to prove.
\end{proof}
}
\noindent With this lemma in hand, we can handle our main sum. We get that 
\begin{equation}\label{eq11}
    S_{\f(\chi)\,\le\, \log^{A} m}=\sum_{\substack{d\le \log^{A} m\\ d\mid m}}\,\,\sum_{\substack{\chi\neq \chi_0,\psi\\\f(\chi)=d}}|S_1(\chi)||S_2(\chi)||S_3(\chi)|.
\end{equation}
By Lemma \ref{lm} and $(\ref{eq6})-(\ref{eq9})$, the terms with non-squarefree $d$ are negligible. {Also, we can observe that $d\ge 5$ is equivalent to the fact that $\chi\neq \chi_{0},\psi$.} Thus, by $(\ref{eq6})-(\ref{eq9})$, we have
\begin{align*}
     S_{\f(\chi)\,\le\, \log^{A} m}&=\sum_{\substack{5\,\le \,d\le \log^{A}m\\ \mu^{2}(d)=1\,,d\mid m}}\,\,\sum_{\substack{\f(\chi)=d}}|S_1(\chi)||S_2(\chi)||S_3(\chi)|+\mathcal{O}_{B}\left(\frac{m^{2+2/k}}{\log ^{B} m}\right)\\&=\sum_{\substack{5\,\le \,d\le \log^{A}m\\ \mu^{2}(d)=1\,,d\mid m}}\,\,\sum_{\substack{\f(\chi)=d}}|I_1||I_2||I_3|\prod_{p\mid d}\left(1+\frac1{p-2}\right)^{3}(|\rho_{\chi_d}|)^{3}+\mathcal{O}_{B}\left(\frac{m^{2+2/k}}{\log^{B-A}m}\right)\\&=|I_1||I_2||I_3|\sum_{\substack{5\,\le \,d\le \log^{A}m\\ \mu^{2}(d)=1\,,d\mid m}}\,\,\sum_{\substack{\f(\chi)=d}}\,\prod_{p\mid d}(p-2)^{-3}+\mathcal{O}_{B}\left(\frac{m^{2+2/k}}{\log^{B-A}m}\right)
\end{align*}

We know that for a square-free integer $d$, there are $\prod_{p\mid d}(p-2)$ primitive characters$\Mod{d}$. We get 
\begin{align}\label{eq12}
S_{\f(\chi)\le\, \log^{A} m}\le |I_1||I_2||I_3|\sum_{\substack{d=5,\,2\nmid d\\\mu^{2}(d)=1}}^{\infty}\,\,\,\prod_{p\mid d}(p-2)^{-2}+\mathcal{O}_{B}\left(\frac{m^{2+2/k}}{\log^{B-A}m}\right), \end{align} where the main term can be bounded using \begin{align*}\sum_{\substack{d=5,\,2\nmid d\\\mu^{2}(d)=1}}^{\infty}\,\,\,\prod_{p\mid d}(p-2)^{-2}&= \prod_{p\ge3}\left(1+\frac{1}{(p-2)^{2}}\right)-2\le 2\exp\left(\sum_{p\ge 5}\frac1{(p-2)^{2}}\right)-2\\&\le 2\exp\left(\sum_{n\ge2}\frac{1}{(2n-1)^{2}}\right)-2\\&\le2\exp(\pi^{2}/8-1)-2=0.526\ldots\end{align*}Finally, we get that 
\begin{equation}\label{eq13}
    S_{\f(\chi)\,\le\, \log^{A} m} \le (0.53+o(1))|I_1||I_2||I_3|,
\end{equation}
once we choose $B=A+4$.
\section{Sum involving characters of large conductor and final steps}
To estimate $S_{\f(\chi)\,>\,\log^{A}m}$, we essentially follow the approach in \cite{Cilleruelo2009}. We will use the lemma below.
\begin{lem}\label{lm3}
    If $\chi\neq \chi_0$ such that $\f(\chi)>\log^{4(k+3)^{2}} m$, then \[S_1(\chi)\ll \frac{|I_1|\log\log m }{\log^{k^{2}+6k+3} m}.\]
\end{lem}
\begin{proof}
    By Rakhmonov's estimate~\cite[Theorem 1]{zbMATH01173130}, we know that for any nontrivial character $\chi\Mod{m}$,
    \[\bigg|\sum_{p\le x}\chi(p-1)\bigg|\le x(\log x)^{5}\tau(q)\left(\sqrt{1/q+q\tau^{2}(q_1)/x}+x^{-1/6}\tau(q_1)\right),\] where $q=\f(\chi)$ and $q_1=\prod_{p\mid m, p\nmid q}p$. 
    We know that \[
|S_1(\chi)| \leq \left| \sum_{p \leq m^{1+1/k}} \chi(p - 1) \right| + \left| \sum_{p \leq 0.5m^{1+1/k}} \chi(p - 1) \right|.
\]

For $x = m^{1+1/k}$ or $x = 0.5m^{1+1/k}$, this gives
\begin{align*}
\sum_{p \leq x} \chi(p - l) &\ll m^{1+1/k} (\log m)^5 \tau(q) / \sqrt{q}  \\& +m^{1/2+1/(2k)} (\log m)^5 q^{1/2} \tau(q_1) \tau(q) \\&+ m^{(1+1/k)5/6} (\log m)^5 \tau(q_1) \tau(q).
\end{align*}
Since $q \leq m$, $k \geq 2$ and $\tau(q_1) \tau(q) \leq \tau(m) \ll m^{1/(4k)}$, we obtain
\[
\sum_{p \leq x} \chi(p - l) \ll m^{1+1/k} (\log m)^5 \tau(q) / \sqrt{q} + m^{1+3/(4k)} (\log m)^5.
\]
We know that $q>\log^{4(k+3)^{2}} m$. Therefore, we get
\begin{align*}
\sum_{p \leq x} \chi(p - l) &\ll m^{1+1/k} (\log m)^{5-(k+3)^2} + m^{1+3/(4k)} (\log m)^5 \\&\ll \frac{m}{ \varphi(m)}{(\log m)^{6-(k+3)^2}} |I_1|.
\end{align*}
Finally, we use the known estimate $\varphi(m) \gg m / \log \log m$.
\end{proof}

\begin{remark}
    Kerr \cite[Theorem 1]{MR4324085} has proven that for a primitive Dirichlet character $\chi\Mod{q}$ and integer $a$ coprime with $q$, 
    \[\bigg| \sum_{n\le N}\Lambda(n)\chi(n+a)\bigg|\,\le\, q^{1/9+o(1)}N^{23/27}\] for $N\ge q$ where $\Lambda(n)$ is the von-Mangoldt function. This result is stronger than the previous Rakhmonov estimate for $N>q^{3/4+\epsilon}$. However, as pointed out in \cite[Section 4]{Cilleruelo2009}, it does not affect the strength of our results.\end{remark}

\noindent The next lemma gives an upper bound on sums involving $S_{j}(\chi)$.
\begin{lem}\label{lm4}
    The following bounds hold:
    \begin{align*}
        &\sum_{\chi}|S_j(\chi)|^{2}\ll |I_j|^{2}\,\log m,\,j=1,2,\\&\sum_{\chi}|S_3(\chi)|^{2k}\ll m^{2}\log^{k^{2}-1} m.
    \end{align*}
\end{lem}
\begin{proof}
    The proof follows exactly as shown in \cite[Lemma 2.2]{Cilleruelo2009}.
\end{proof}

\noindent Employing the Lemmas \ref{lm3} and \ref{lm4} along with H\"older's inequality, we have 
\begin{align*}
S_{\f(\chi)\,>\,\log^{A}m}&\ll \bigg(\max_{\f(\chi)>\log^{A}m}|S_{1}(\chi)|\bigg)^{1/k} \sum_{\substack{\f(\chi)>\log^{A}m}}|S_1(\chi)|^{1-1/k}\,|S_2(\chi)|\,|S_3(\chi)|\\& \ll \frac{|I_1|^{1/k} }{\log^{k+6} m} \,g_1(k)\,g_2(k)\,g_3(k).
\end{align*}
\newline
Here, the sum $g_i(k)$ is defined as following:
\begin{align*}
    &g_1(k)=\big(\sum_{\chi}|S_1(\chi)|^2\big)^{1/2-1/(2k)} \ll |I_1|^{1-1/k}\,(\log m)^{1/2},\\&g_2(k)=\big(\sum_{\chi}|S_2(\chi)|^2\big)^{1/2} \ll |I_2|\,(\log m)^{1/2},\\&g_3(k)=\big(\sum_{\chi}|S_1(\chi)|^{2k}\big)^{1/2k} \ll m^{1/k}\,(\log m)^{k/2}.
\end{align*}

\noindent Thus, we get that  \begin{equation}\label{eq15}
S_{\f(\chi)\,>\,\log^{A}m} \ll \frac{|I_1||I_2|m^{1/k}}{\log^{k/2+5}m} \ll \frac{|I_1||I_2||I_3|}{\log^{k/2+2} m}.
\end{equation}
Combining the estimates $(\ref{eq3})$, $(\ref{eq13})$ and $(\ref{eq15})$, we can see that \begin{align*}
S=\sum_{\chi \neq \chi_0,\psi} |S_1(\chi)||S_2(\chi)||S_3(\chi)| &< (0.53+o(1))|I_1||I_2||I_3| + \mathcal{O}\left(\frac{|I_1||I_2||I_3|}{\log^{k/2+2} m}\right)\\& < |I_1||I_2||I_3|
\end{align*}
for large $m$.

\bibliographystyle{amsplain-no-bysame}
\bibliography{main}

\providecommand{\bysame}{\leavevmode\hbox to3em{\hrulefill}\thinspace}
\providecommand{\MR}{\relax\ifhmode\unskip\space\fi MR }
\providecommand{\MRhref}[2]{%
  \href{http://www.ams.org/mathscinet-getitem?mr=#1}{#2}
}
\providecommand{\href}[2]{#2}
\begin{thebibliography}{10}

\bibitem{Cilleruelo2009}
Javier Cilleruelo and Moubariz~Z. Garaev, \emph{Least totients in arithmetic progressions}, Proc. Amer. Math. Soc. \textbf{137} (2009), no.~9, 2913--2919.

\bibitem{MR1642868}
Thomas Dence and Carl Pomerance, \emph{Euler's function in residue classes}, Ramanujan J. \textbf{2} (1998), no.~1-2, 7--20.

\bibitem{MR1689545}
Kevin Ford, Sergei Konyagin, and Carl Pomerance, \emph{Residue classes free of values of {E}uler's function}, Number theory in progress, {V}ol. 2 ({Z}akopane-{K}o\'scielisko, 1997), de Gruyter, Berlin, 1999, pp.~805--812.

\bibitem{MR2471950}
John~B. Friedlander and Florian Luca, \emph{Residue classes having tardy totients}, Bull. Lond. Math. Soc. \textbf{40} (2008), no.~6, 1007--1016.

\bibitem{MR2331570}
John~B. Friedlander and Igor~E. Shparlinski, \emph{Least totient in a residue class}, Bull. Lond. Math. Soc. \textbf{39} (2007), no.~3, 425--432.

\bibitem{MR2418809}
John~B. Friedlander and Igor~E. Shparlinski, \emph{Corrigendum: ``{L}east totient in a residue class'' [{B}ull. {L}ond. {M}ath. {S}oc. {\bf 39} (2007), no. 3, 425--432 ]}, Bull. Lond. Math. Soc. \textbf{40} (2008), no.~3, 532.

\bibitem{MR2506368}
Moubariz~Z. Garaev, \emph{A note on the least totient of a residue class}, Q. J. Math. \textbf{60} (2009), no.~1, 53--56.

\bibitem{MR4324085}
Bryce Kerr, \emph{Bounds of multiplicative character sums over shifted primes}, Tr. Mat. Inst. Steklova \textbf{314} (2021), 71--96.

\bibitem{MR3971232}
Dimitris Koukoulopoulos, \emph{The distribution of prime numbers}, Graduate Studies in Mathematics, vol. 203, American Mathematical Society, Providence, RI, [2019] \copyright 2019.

\bibitem{MR3356031}
Youness Lamzouri, Xiannan Li, and Kannan Soundararajan, \emph{Conditional bounds for the least quadratic non-residue and related problems}, Math. Comp. \textbf{84} (2015), no.~295, 2391--2412.

\bibitem{pollack2024mean}
Paul Pollack and Akash~Singha Roy, \emph{Mean values of multiplicative functions and applications to residue-class distribution}, Proc. Edinburgh Math. Soc. (2024), 1--19.

\bibitem{zbMATH01173130}
Z.~Kh. Rakhmonov, \emph{On the distribution of the values of {D}irichlet characters and their applications}, Trudy Mat. Inst. Steklov. \textbf{207} (1994), 286--296.

\bibitem{MR2825574}
Triantafyllos Xylouris, \emph{On the least prime in an arithmetic progression and estimates for the zeros of {D}irichlet {$L$}-functions}, Acta Arith. \textbf{150} (2011), no.~1, 65--91.

\end{thebibliography}

\end{document}